\documentclass{article}
 \usepackage{amsmath}

\newtheorem{thm}{Theorem}[section]
\newtheorem{cor}[thm]{Corollary}
\newtheorem{prop}[thm]{Proposition}

\newtheorem{lem}[thm]{Lemma}

\newtheorem{rem}[thm]{Remark}

\newtheorem{ex}[thm]{Example}

\newcommand{\be}{\begin{equation}}
\newcommand{\ee}{\end{equation}}
\newcommand{\ben}{\begin{enumerate}}
\newcommand{\een}{\end{enumerate}}
\newcommand{\beq}{\begin{eqnarray}}
\newcommand{\eeq}{\end{eqnarray}}
\newcommand{\beqn}{\begin{eqnarray*}}
\newcommand{\eeqn}{\end{eqnarray*}}

\newcommand{\pa}{\partial}

\newcommand{\pxi}{ {\pa \over \pa x^i}}

\newcommand{\pyj}{ {\pa \over \pa y^j}}

\newcommand{\qed}{\hspace*{\fill}Q.E.D.}  

\begin{document}
\title{Conformal Vector Fields On Projectively Flat $(\alpha,\beta)$-Finsler Spaces }
\author{Guojun Yang\footnote{Supported by the
National Natural Science Foundation of China (11471226) }}
\date{}
\maketitle
\begin{abstract}
 In this paper, it is proved that any
conformal vector field is homothetic on a locally projectively
flat $(\alpha,\beta)$-space of non-Randers type in dimension $n\ge
3$, and the local solutions of such a vector field are determined.
While on a locally projectively flat Randers space, examples show
that the conformal vector fields are not necessarily homothetic.


 {\bf MR(2000) subject classification: }
53B40, 53C60
\end{abstract}

\section{Introduction}

Conformal vector fields play an important role in Finsler
geometry.  When $F$ is a Riemann metric, the local solutions of a
conformal vector field can be determined
 if $F$ satisfies certain curvature conditions (cf. \cite{YShen} \cite{ShX}
 \cite{Y1} \cite{Y2}). Some  problems in Finsler geometry can be
solved by constructing a conformal vector field of a Riemann
metric with certain curvature features (cf. \cite{BRS}
\cite{Y3}--\cite{Yu1}). For two conformally related Finsler
metrics on a manifold, their conformal vector fields coincide
(Proposition \ref{prop022} below).

A conformal vector field has some different equivalent definitions
(see Section \ref{sec2} below). Every conformal vector field $V$
is associated with a scalar function $c$ called the conformal
factor. If $c$ is a constant, $V$ is said to be homothetic; if
$c=0$, $V$ is said to be Killing. As a special case of conformal
vector fields, homothetic vector fields have some special
properties. For example, Huang-Mo obtain the relation between the
flag curvatures of two Finsler metrics $F$ and $\tilde{F}$, where
$\tilde{F}$ is defined by $(F,V)$ under navigation technique for a
homothetic vector field $V$ of $F$ (\cite{MH}).

A Randers metric is defined by $F=\alpha+\beta$, where $\alpha$ is
a Riemann metric and $\beta$ is a 1-form. In \cite{SX}, Shen-Xia
study conformal vector fields of Randers spaces under certain
curvature conditions. In \cite{HM1}, Huang-Mo show that a
conformal vector field of a Randers space of isotropic S-curvature
must be homothetic. As a generalization of Randers metrics,
 an  $(\alpha,\beta)$-metric is defined by
 $$F=\alpha \phi(s),\ \ s=\beta/\alpha,$$
where $\alpha$ is a  Riemann metric,  $\beta$ is a 1-form and
$\phi(s)$ is a function satisfying certain conditions. An
$(\alpha,\beta)$-metric $F=\alpha\phi(\beta/\alpha)$ is said to be
of Randers type, if $\phi(s)= \sqrt{1+ks^2}+\epsilon s$ for some
constants $k,\epsilon$.   In \cite{K}, Kang characterizes the
conformal vector fields of an $(\alpha,\beta)$-space  by some PDEs
in a special case $\phi'(0)\ne 0$. Actually, we can obtain the
same result as in \cite{K} without the condition $\phi'(0)\ne 0$
and it only requires that the $(\alpha,\beta)$-space be
non-Riemannian (see Proposition \ref{th1} below).

It is the Hilbert's Fourth Problem to study   projectively flat
metrics. A Finsler metric is said to be locally projectively flat
if its geodesics are straight in suitable local coordinate
systems. A Randers metric $F=\alpha+\beta$ is locally projectively
flat if and only if $\alpha$ is of constant sectional curvature
and $\beta$ is closed (\cite{BaMa}). Shen characterizes
projectively flat $(\alpha,\beta)$-metrics of non-Randers type
(\cite{Shen1}). In this paper, we will study conformal vector
fields of projectively flat $(\alpha,\beta)$-spaces.

Let  $F=\alpha\phi(\beta/\alpha)$ be an $n(\ge 3)$-dimensional
locally projectively flat $(\alpha,\beta)$-metric of non-Randers
type. It is known that $F$ can be determined by a Riemann metric
$h$ of constant sectional curvature $\mu$
  and a closed and conformal 1-form $\rho$ with respect to $h$ which
  can be
locally expressed as (Lemma \ref{lem23} below, or cf. \cite{Shen1}
\cite{Y7} \cite{Yu1} \cite{Y1})
 \be\label{y2}
 h=\frac{2}{1+\mu|x|^2}|y|, \ \ \ \ \
 \rho=\frac{4}{(1+\mu|x|^2)^2}\Big\{-2(\lambda+\mu \langle
 e,x\rangle)\langle x,y\rangle+(1+\mu|x|^2)\langle
 e,y\rangle\Big\},
 \ee
where $\lambda,\mu$ are constant numbers, and $e=(e^i)$ is a
constant vector. Note that when $F$ is given, $\lambda,\mu,e$
 in (\ref{y2}) are determined. We will use $\lambda,\mu,e$ in Theorem \ref{th2}
below.

\begin{thm}\label{th2}
  Let  $F=\alpha\phi(\beta/\alpha)$ be an $n(\ge 3)$-dimensional
locally projectively flat $(\alpha,\beta)$-metric of non-Randers
type with $\phi(0)=1$.  Then any conformal vector field $V$ of $F$
is homothetic. Further, put $V=V^i(x)\pa/\pa x^i$, and then
related to (\ref{y2}),
 locally  we have one of the following  cases:
 \ben
  \item[{\rm (i)}] ($\mu=0$, $\lambda=0$) $V$ is given by
   \be\label{y3}
 V^i=-2\tau x^i+q^i_kx^k+\gamma^i,\ \ \ (Qe=0),
   \ee
  In this case,   $F$ is
  flat-parallel and
  $c$ in (\ref{y1}) is given by $c=\tau$.

\item[{\rm (ii)}] ($\mu=0$, $\lambda\ne 0$) $V$ is given by
  \be\label{y03}
 V^i=q^i_kx^k+\gamma^i,\ \ \ (Qe=-2\lambda \gamma),
   \ee
 In this case,
  $c$ in (\ref{y1}) is given by $c=0$, and so $V$ is a Killing
  field.

 \item[{\rm (iii)}] ($\mu\ne 0$)  $V$ is given by
   \be\label{y4}
 V^i=2\mu\langle \gamma,x\rangle
 x^i+(1-\mu|x|^2)\gamma^i+q^i_kx^k,\ \ \ \big(\langle\gamma,e\rangle=0,\  Qe=-2\lambda \gamma\big),
   \ee
  In this case,
  $c$ in (\ref{y1}) is given by $c=0$, and so $V$ is a Killing
  field.
 \een
 In the above,  $\tau$ is a constant number, $\gamma=(\gamma^i)$ is a constant
 vector, and $Q=(q^i_k)$ is a constant skew-symmetric matrix.
\end{thm}

Note that for $\lambda\ne 0$ in (\ref{y4}), we have
$Qe=-2\lambda\gamma  \Longrightarrow \langle\gamma,e\rangle=0$.
 Theorem \ref{th2} is not true if
$F$ is a metric of Randers type. For a locally projectively flat
Randers space, its conformal vector field   is not necessarily
homothetic. The following is an example.

\begin{ex}
 Let $F=\alpha+\beta$ be a Randers metric on a manifold $M$. Locally, define
  \beqn
 \alpha:&&\hspace{-0.6cm}=\frac{2}{1+\mu|x|^2}|y|, \ \
 \beta:=\frac{1}{\lambda(1-\mu|x|^2)+\langle
 d,x\rangle}\Big\{\langle d,y\rangle-\frac{2\mu(2\lambda+\langle d,x\rangle)\langle x,
 y\rangle}{1+\mu|x|^2}\Big\},\\
 V^i:&&\hspace{-0.6cm}=-2\big(\lambda+\langle
 d,x\rangle\big)x^i+|x|^2d^i,
  \eeqn
  where the constant $\lambda$ and the constant vector $d=(d^i)\ne 0$
  satisfy $|d|^2+4\mu\lambda^2=0$. It is easy to
  see that $F$ is locally projectively flat.

  By Lemma \ref{lem22} (i)
  below, the vector field $V=(V^i)$ is conformal on $(M,\alpha)$.
  It can be directly verified that $\beta$ (closed) and $V$ satisfy  (\ref{y1}) with
   $$c=\frac{\lambda(1-\mu|x|^2)+\langle d,x\rangle}{1+\mu|x|^2}\ (\ne
   constant).$$
   Thus by Proposition \ref{th1}, $V$ is a conformal vector field of
   $(M,F)$ with the conformal factor $c$, and $V$ is not homothetc.
\end{ex}

We show in \cite{SY} \cite{Y7} that $F=\alpha+\epsilon \beta\pm
\beta^2/\alpha$ for a constant $\epsilon$  is of scalar flag
curvature in dimension $n\ge 3$ if and only if $F$ is locally
projectively flat. By Theorem \ref{th2}, we have the following
corollary.

\begin{cor}\label{cor3}
 Let $F=\alpha+\epsilon \beta\pm
\beta^2/\alpha$ be an  $n(\ge 3)$-dimensional
$(\alpha,\beta)$-metric of scalar flag curvature, where $\epsilon$
is a constant. Suppose $V=V^i(x)\pa/\pa x^i$ is a
 conformal vector field of $F$. Then $V$ can be locally written in the
 form (\ref{y3}), (\ref{y03}), or (\ref{y4}).
\end{cor}

In Theorem \ref{th2}, the result remains the same for a locally
projectively flat Randers metric of isotropic S-curvature
(Corollary \ref{cor41} below). We also give the local structure of
conformal vector fields of a locally projectively flat
$(\alpha,\beta)$-metric with constant flag curvature (Corollary
\ref{cor42} below).

\section{Conformal vector fields}\label{sec2}

Let $F$ be a Finsler metric on  a manifold $M$, and $V$ be a
vector field on $M$. Let $\varphi_t$ be the flow generated by $V$.
Define $\widetilde{\varphi}_t: TM \mapsto TM$ by
  $\widetilde{\varphi}_t(x,y)=(\varphi_t(x),\varphi_{t*}(y))$.
  $V$ is said to be conformal if (cf. \cite{HM1})
 \be\label{j8}
 \widetilde{\varphi}_t^*F=e^{-2\sigma_t}F,
 \ee
where $\sigma_t$ is a function on $M$ for every $t$.
 Differentiating (\ref{j8}) by $t$ at $t=0$,  we obtain
  \be\label{y6}
 X_V(F)=-2cF,
  \ee
  where we define
   \be\label{XV}
   X_V:=V^i\pxi+y^i\frac{\pa
   V^j}{\pa x^i}\pyj,\ \ \ \ \  c:=\frac{d}{dt}|_{t=0}\sigma_t.
   \ee
 In (\ref{XV}), the function $c$ is called the conformal factor.

   \begin{rem}\label{rem1}
(\cite{HM1}) A vector field $V$ is conformal satisfying (\ref{j8})
if and only if (\ref{y6}) holds for some scalar function $c$. In
this
 case, $c$ and $\sigma_t$ are related by
 \be\label{gg9}
 \sigma_t=\int_0^t c(\varphi_s)ds, \ \ \ c=\frac{d}{dt}_{|t=0}\sigma_t.
 \ee
\end{rem}

\begin{rem}\label{rem2}
By (\ref{gg9}), we easily see that
$\widetilde{\varphi}_t^*F=e^{-2ct}F$ for a scalar function $c$ if
and only if $c$ is constant along every integral curve of $V$.

\end{rem}

\begin{lem}\label{lem20}
A vector
   field $V$ on a Finsler manifold $(M,F)$ is conformal with the conformal factor $c$  if and
   only if
    \be\label{y011}
  X_V(F^2)=-4cF^2  \ (\Longleftrightarrow
 X_V(F)=-2cF),\ \ \  or \  V_{0|0}=-2cF^2.
   \ee
   where $_|$ is the $h$-covariant derivative  of Cratan (Berwald, or Chern) connection.
\end{lem}

{\it Proof :} By Remark \ref{rem1}, we only need to prove
 $X_V(F^2)=2V_{0|0}$.
Let
 $$
 G_i=g_{im}G^m=\frac{1}{4}\big
 \{[F^2]_{x^ky^i}y^k-[F^2]_{x^i}\big\},\ \ G^k_i=\frac{\pa G^k}{\pa
 y^i},\ \  G^k_{ij}=\frac{\pa^2 G^k}{\pa y^i\pa y^j}.
 $$
Then a direct computation shows that
 $$
X_V(F^2)=V^i(F^2)_{x^i}+2\frac{\pa V^i}{\pa
x^j}y_iy^j=V^i(F^2)_{x^i}+2\frac{\pa V_i}{\pa
x^j}y^iy^j-2\frac{\pa y_i}{\pa x^j}V^iy^j,
 $$
 $$
\frac{\pa V_i}{\pa x^j}y^iy^j=V_{0|0}+(G^k_j\frac{\pa V_i}{\pa
y^k}+V_kG^k_{ij})y^iy^j=V_{0|0}+2G^k\frac{\pa V_i}{\pa
y^k}y^i+2V_kG^k,
 $$
 $$
 4V_kG^k=4V^kG_k=V^i\{[F^2]_{x^ky^i}y^k-[F^2]_{x^i}\big\}=2V^i\frac{\pa y_i}{\pa
 x^k}y^k-V^i(F^2)_{x^i}.
 $$
In the above three equations, plugging the second and third
equations into the first one gives
 $$
X_V(F^2)=2V_{0|0}+4G^k\frac{\pa V_i}{\pa y^k}y^i=2V_{0|0},
 $$
where we have used $\frac{\pa V_i}{\pa y^k}y^i=0$.   \qed

\begin{lem}\label{lem21}
Let $\beta=b_i(x)y^i$ be a 1-form, and $V$ be a vector field on a
Riemann manifold $(M,\alpha)$ with $\alpha=\sqrt{a_{ij}y^iy^j}$.
Then we have
 \be\label{y8}
 X_V(\alpha^2)=2V_{0;0},\ \ \ \ \ X_V(\beta)=(V^j\frac{\pa b_i}{\pa
 x^j}+b_j\frac{\pa V^j}{\pa x^i})y^i=(V^jb_{i;j}+b^jV_{j;i})y^i,
 \ee
 where $V_i:=a_{ij}V^j$ and $b^i:=a^{ij}b_j$, and the covariant derivative is taken with respect to the
 Levi-Civita connection of $\alpha$.
\end{lem}

\begin{lem} \label{lem22}
Let $\alpha$ be a Riemann metric of constant sectional curvature
$\mu$ on an $n$-dimensional manifold $M$. Locally express $\alpha$
by
 \be\label{g14}
\alpha=\frac{2}{1+\mu|x|^2}|y|.
 \ee
 \ben
  \item[{\rm (i)}] (\cite{Y1}) ($n\ge 3$) Let $V$ be a
conformal vector field of $(M,\alpha)$ with the conformal factor
$c=c(x)$. Then locally we have
 \be\label{g15}
 V^i=-2\big(\lambda+\langle
 d,x\rangle\big)x^i+|x|^2d^i+q_r^ix^r+\eta^i,\ \
c=\frac{\lambda(1-\mu |x|^2)+ \langle \mu \eta+d,x\rangle}{1+\mu
 |x|^2},
 \ee
where $\lambda$ is a constant number, $d,\eta$ are constant
vectors and
  $(q_i^j)$ is skew-symmetric.
\item[{\rm (ii)}] ($n\ge 2$) In (i), if additionally the 1-form
$V_iy^i$ is closed ($V_i:=a_{ik}V^k$), then
 \be\label{g16}
 V^i=-2(\lambda+\mu \langle
 e,x\rangle) x^i+(1+\mu|x|^2)
 e^i,\ \ \ c=\frac{\lambda(1-\mu |x|^2)+ 2\mu\langle e,x\rangle}{1+\mu
 |x|^2}.
 \ee
 \een
\end{lem}

{\it Proof :}  We prove in \cite{Y1} that if $\alpha$ is locally
conformally flat taking the local form
 \be\label{g17}
 \alpha=e^{\frac{1}{2}\sigma(x)}|y|,
 \ee
 then $V$ is a conformal vector field of $\alpha$ if and only if
 \be\label{g18}
 \frac{\pa V^i}{\pa x^j}+\frac{\pa V^j}{\pa
  x^i}=0 \quad{\rm (} \forall \thinspace i\ne  j{\rm )},\ \ \  \frac{\pa V^i}{\pa x^i}=\frac{\pa V^j}{\pa
  x^j} \quad{\rm (} \forall \thinspace i, j{\rm )}.
 \ee
In this case, $c$ is given by
 \be\label{g19}
 c(x)=-\frac{1}{4}\big[2\frac{\pa V^1}{\pa x^1}+V^r\sigma_r\big],\
 \ \sigma_i:=\sigma_{x^i}.
 \ee
In \cite{Y1}, we show that when $n\ge3$, the solutions $V^i$ of
(\ref{g18}) are given by (\ref{g15}). If  $\alpha$ takes the local
form
  (\ref{g14}), then  $\sigma$ in (\ref{g17}) is given by
   $$
 \sigma=ln\frac{4}{(1+\mu|x|^2)^2}.
   $$

For Case (i),  by (\ref{g19}) and $V$ in (\ref{g15}) we get $c$
given by (\ref{g15}).

For Case (ii), if $n\ge 3$, then since $V_0$ is closed, it is easy
to show that $Q=0$ and $d=\mu \eta$ in (\ref{g15}). Thus we get
$V$ and $c$ given by (\ref{g16}) (put $\eta=e$) from (\ref{g15}).
If $n=2$, then since $\pa V^1/\pa x^2=-\pa V^2/\pa x^1$ from
(\ref{g18}), it follows that $V_0$ is closed if and only if
 \be\label{g20}
\frac{\pa V^2}{\pa x^1}=\frac{2\mu(V^2x^1-V^1x^2)}{1+\mu|x|^2}.
 \ee
Solving the system (\ref{g18}) and (\ref{g20}) yields $V$ in
(\ref{g16}). Then it follows from (\ref{g19}) that $c$ is given by
(\ref{g16}). \qed

\begin{prop}\label{prop022}
   Let $(M,F)$ and
  $(M,\widetilde{F})$ be conformally related with
  $\widetilde{F}=e^{\sigma/2}F$ for a scalar function $\sigma$. Then
 $V$ is a conformal vector field of  $(M,F)$ iff. $V$ is a conformal vector field of
 $(M,\widetilde{F})$. Further, their conformal factors $c$ and
 $\widetilde{c}$ are related by
  \be\label{cc}
  \widetilde{c}=c-\frac{1}{4}\ V(\sigma).
  \ee
\end{prop}

{\it Proof :} Assume that $V$ is a conformal vector field of
$(M,F)$ with the conformal factor $c$. Then by Lemma \ref{lem20}
we have
  $$
 X_V(F^2)=-4cF^2.
  $$
Thus we have
 $$
 X_V(\widetilde{F}^2)=X_V(e^{\sigma}F^2)=X_V(e^{\sigma})F^2+e^{\sigma}X_V(F^2)=\big[V(\sigma)-4c\big]\widetilde{F}^2,
 $$
which,  by Lemma \ref{lem20} again, implies that $V$ is also a
conformal vector field of $(M,\widetilde{F})$, and the conformal
factor $\widetilde{c}$ is given by (\ref{cc}).  \qed

\

A conformal vector field of an $(\alpha,\beta)$-metric is
characterized in the following proposition, which has been proved
by Kang in \cite{K} under an additional condition $\phi'(0)\ne0$.

\begin{prop}\label{th1}
  Let $F=\alpha\phi(\beta/\alpha)$ be a non-Riemannian
  $(\alpha,\beta)$-metric with
 $\phi(0)=1$, and $V=V^i(x)\pa/\pa x^i$ be a vector field.
 Then $V$ is a conformal vector field of $F$ with the conformal factor $c$ if and only if
  \be\label{y1}
  V_{i;j}+V_{j:i}=-4ca_{ij},\ \ \ \ V^jb_{i;j}+b^jV_{j;i}=-2cb_i,
  \ee
where $c$ is a scalar function, $V_i$ and $b^i$ are defined by
$V_i:=a_{ij}V^j$ and $b^i:=a^{ij}b_j$, and the covariant
derivatives are taken with respect to the Levi-Civita connection
of $\alpha$.
\end{prop}

{\it Proof :} By Lemma \ref{lem20} ,  $V$ is a conformal vector
field of $F$ iff. $X_V(F^2)=-4cF^2$. A direct computation shows
that
 \beqn
 X_V(F^2)&=&\phi^2X_V(\alpha^2)+2\alpha^2\phi\phi'\frac{\alpha
 X_V(\beta)-\beta
 X_V(\alpha)}{\alpha^2}\\
 &=&\phi(\phi-s\phi')X_V(\alpha^2)+2\alpha\phi\phi'X_V(\beta).
 \eeqn
Now plugging (\ref{y8}) into the above equation, we see that
$X_V(F^2)=-4cF^2$ is written as
 \be\label{y20}
 V_{0;0}+\alpha
 Q(V^ib_{j;i}+b^iV_{i;j})y^j=-\frac{2c\phi}{\phi-s\phi'}\alpha^2,\
 \ \ \  (Q:=\frac{\phi'}{\phi-s\phi'}),
  \ee
where the covariant derivative is taken with respect to
  $\alpha$.

In order to simplify (\ref{y20}), we choose a special coordinate
 system $(s,y^a)$ at  a fixed point on a manifold as usually used.
Fix an arbitrary point $x\in M$ and take  an orthogonal basis
  $\{e_i\}$ at $x$ such that
   $$\alpha=\sqrt{\sum_{i=1}^n(y^i)^2},\ \ \beta=by^1.$$
It follows from $\beta=s\alpha$ that
 $$y^1=\frac{s}{\sqrt{b^2-s^2}}\bar{\alpha},\ \ \ \ \Big(\bar{\alpha}:=\sqrt{\sum_{a=2}^n(y^a)^2}\Big).$$
Then if we  change coordinates $(y^i)$ to $(s, y^a)$, we get
  $$\alpha=\frac{b}{\sqrt{b^2-s^2}}\bar{\alpha},\ \
  \beta=\frac{bs}{\sqrt{b^2-s^2}}\bar{\alpha}.$$
 Let
 $$\bar{V}_{0;0}:=V_{a;b}y^ay^b, \ \ \bar{V}_{1;0}:=V_{1;a}y^a,\ \ \bar{V}_{0;1}:=V_{a;1}y^a,
 \ \ \bar{b}_{0;i}:=b_{a;i}y^a $$
Note that under the coordinate $(s,y^a)$, we have
$b_1=b,\bar{b}_0=0$, but generally $\bar{b}_{0;i}\ne 0$.

Under the coordinate $(s,y^a)$, (\ref{y20}) is equivalent to
  \beq
 0&=&b(b\bar{V}_{1;0}+V^i\bar{b}_{0;i})Q+
 (\bar{V}_{1;0}+\bar{V}_{0;1})s,\label{y23}\\
 0&=&\big[b(2bc+V^ib_{1|i}+bV_{1|1})sQ+2b^2c+V_{1|1}s^2\big]\bar{\alpha}^2+(b^2-s^2)\bar{V}_{0|0}.\label{y24}
  \eeq

For (\ref{y23}), we will prove
 \be\label{yy0023}
b\bar{V}_{1;0}+V^i\bar{b}_{0;i}=0,\ \ \ \
 \bar{V}_{1;0}+\bar{V}_{0;1}=0.
 \ee
If $b\bar{V}_{1;0}+V^i\bar{b}_{0;i}\ne 0$ at a point, then we have
$Q=ks$ for some constant $k$. Solving $Q=ks$ with $\phi(0)=1$ we
get $\phi(s)=\sqrt{1+ks}$ and hence $F$ is Riemannian. So we must
have the first equation in (\ref{yy0023}). By (\ref{y23}) again,
we have the second equation in (\ref{yy0023}).

For (\ref{y24}), we will prove
 \be\label{yy0024}
\bar{V}_{0;0}=-2c\bar{\alpha}^2,\ \ \  V^ib_{1|i}+bV_{1|1}=-2bc,\
\ \  V_{1|1}=-2c.
 \ee
Putting $s=0$ in (\ref{y24}) we have
$\bar{V}_{0;0}=-2c\bar{\alpha}^2$. By this, (\ref{y24}) is reduced
to
 \be\label{yy0025}
 b(2bc+V^ib_{1|i}+bV_{1|1})Q+(V_{1|1}+2c)s=0.
 \ee
Similar to the analysis of (\ref{y23}), we  obtain the second and
third equations of (\ref{yy0024}).

Finally, (\ref{y1}) immediately follows from (\ref{yy0023}) and
(\ref{yy0024}).          \qed

\section{Proof of Theorem \ref{th2}}\label{sec4}

To prove Theorem \ref{th2}, we first show the following  key
Lemmas \ref{lem23} and \ref{lem41}.

\begin{lem}
(\cite{Shen1} \cite{Y7} \cite{Yu1})\label{lem23}
 Let $F=\alpha \phi(s)$, $s=\beta/\alpha$, be an $n(\ge 3)$-dimensional
  $(\alpha,\beta)$-metric, where $\phi(0)=1$. Suppose that
    $\beta$ is not parallel with respect to $\alpha$ and $F$ is not
of Randers type. Then $F$ is  locally projectively flat if and
only if $\phi(s)$ satisfies
   $$
    \big\{1+(k_1+k_3)s^2+k_2s^4\big\}\phi''(s)=(k_1+k_2s^2)\big\{\phi(s)-s\phi'(s)\big\},
    $$
 and $\alpha,\beta$ are determined by
  \be\label{y02}
 h=\sqrt{u(b^2)\alpha^2+v(b^2)\beta^2},\ \ \rho=w(b^2)\beta,\ \
 (b^2:=||\beta||_{\alpha}^2),
 \ee
   where  $k_1,k_2,k_3$ are constants satisfying $k_2\ne
   k_1k_3$,
   $h$ is a Riemann metric of constant sectional
curvature and $\rho$ is a closed conformal 1-form with respect to
$h$, and $u=u(t)\ne 0,v=v(t),w=w(t)\ne 0$ satisfy the following
ODEs:
 \beq
u'&=&\frac{v-k_1u}{1+(k_1+k_3)t+k_2t^2},\label{u}\\
v'&=&\frac{u(k_2u-k_3v-2k_1v)+2v^2}{u[1+(k_1+k_3)t+k_2t^2]},\label{v}\\
w'&=&\frac{w(3v-k_3u-2k_1u)}{2u\big[1+(k_1+k_3)t+k_2t^2\big]}.\label{w}
 \eeq
\end{lem}

\begin{rem}
We can have different suitable choices of $u,v,w$ satisfying
(\ref{u})--(\ref{w}).
  A suitable  choice of the
  triple $(u,v,w)$ can be taken as (\cite{Yu1})
   $$
   u=e^{2\chi}, \ \ v=(k_1+k_3+k_2b^2)u, \ \ w=\sqrt{1+(k_1+k_3)b^2+k_2b^4}\ e^{\chi},
  $$
   where $\chi$ is defined by
 $$2\chi:=\int_0^{b^2}\frac{k_2t+k_3}{1+(k_1+k_3)t+k_2t^2}dt.$$
\end{rem}

\begin{lem}\label{lem41}
Let $\alpha=\sqrt{a_{ij}y^iy^j}$ be a Riemann metric and
 $\beta=b_iy^i$ be a 1-form and $V=V^i\pa/\pa x^i$ be a vector field
  on an $n$-dimensional manifold $M$.  Define a
 Riemann metric $h=\sqrt{h_{ij}(x)y^iy^j}$ and a 1-form $\rho=p_i(x)y^i$ by
  \be\label{y37}
h=\sqrt{u(b^2)\alpha^2+v(b^2)\beta^2},\ \ \rho=w(b^2)\beta,\ \
 (b^2:=||\beta||_{\alpha}^2),
  \ee
 where $u=u(t)\ne 0,v=v(t),w=w(t)\ne 0$ are arbitrary smooth functions.
  Then $\alpha$, $\beta$  and $V$ satisfy (\ref{y1}) if and
 only if
  \be\label{y38}
 V_{0|0}=-2ch^2,\ \ \ \ V^jp_{i|j}+p^jV_{j|i}=-2cp_i,
  \ee
  where $p^j:=h^{ij}p_i,\ V_j:=h_{ij}V^i$ and the covariant derivative is taken
  with respect to the Levvi-Civita connection of $h$.
\end{lem}

{\it Proof :} Assume (\ref{y1}) holds. Let $X_V$ be defined by
(\ref{XV}). We have
 \beq\label{y40}
 X_V(b^2)&=&2V^jb^ib_{i;j}\nonumber\\
 &=&2(-2cb^2-b^ib^jV_{i;j}),\ \ \text{(by the second equation of
 (\ref{y1}))}\nonumber\\
 &=&0,\ \ \ \ \ \ \ \ \text{(by the first equation of
 (\ref{y1}))}.
 \eeq
By (\ref{y1}) and (\ref{y8}) we get
 \be\label{y41}
 X_V(\alpha^2)=-4c\alpha^2,\ \ \ \ X_V(\beta)=-2c\beta.
 \ee
Therefore, it follows from (\ref{y37}), (\ref{y40}) and
(\ref{y41}) that
 \beq
 X_V(h^2)&=&u'(b^2)X_V(b^2)\alpha^2+v'(b^2)X_V(b^2)\beta^2+u(b^2)X_V(\alpha^2)+v(b^2)X_V(\beta^2)\nonumber\\
 &=&u(b^2)X_V(\alpha^2)+v(b^2)X_V(\beta^2)\nonumber\\
 &=&-4c\big[u(b^2)\alpha^2+v(b^2)\beta^2\big]\nonumber\\
 &=&-4ch^2,\label{y42}\\
X_V(\rho)&=&w'(b^2)X_V(b^2)\beta+w(b^2)X_V(\beta)=w(b^2)X_V(\beta)\nonumber\\
&=&-2cw(b^2)\beta=-2c\rho.\label{y43}
 \eeq
Now (\ref{y38}) follows directly from (\ref{y42}), (\ref{y43}) and
Lemma \ref{lem21}.

Conversely, assume (\ref{y38}) holds. Then we can prove that
(\ref{y1}) holds in a similar way as  above, since we can express
$\alpha$ and $\beta$ in terms of $h$, $\rho$ and the norm of
$\rho$ with respect to $h$. This completes the proof of the lemma.
 \qed

\begin{rem}
 In (\ref{y37}), take $u(b^2)=1-b^2,\  v(b^2)=b^2-1,\
 w(b^2)=b^2-1$,
 and then $(h,\rho)$ is called the navigation data for a Randers
 metric $F=\alpha+\beta$ under navigation representation.
 The method used in the proof of Lemma \ref{lem41} also gives a
 brief proof to Proposition 3.1 in \cite{SX} for the case of
 Randers metrics (cf. the long proof therein).
\end{rem}

{\it Proof of Theorem \ref{th2} :} Let
$F=\alpha\phi(\beta/\alpha)$ be an $(\alpha,\beta)$-metric (not of
Randers type) satisfying the conditions in Theorem \ref{th2}.

\

\noindent{\bf Condition (i):} Suppose $\beta$ is not parallel with
respect to $\alpha$.

 Define $h$ and $\rho$ by (\ref{y37}) (or
(\ref{y02})), where $u=u(t)\ne0,v=v(t),w=w(t)\ne0$ are a suitable
choice satisfying the ODEs (\ref{u})--(\ref{w}). Then by Lemma
\ref{lem23}, $h$ is a Riemann metric of constant sectional
curvature $\mu$ and $\rho$ is a closed conformal 1-form. Since $V$
is a conformal vector field of $F$, $V$ satisfies (\ref{y1}) by
Proposition \ref{th1}. Then by Lemma \ref{lem41}, $V$ satisfies
(\ref{y38}).

 Since $h$ is  of constant sectional
curvature $\mu$ and $\rho$ is  closed and conformal,  $h$ and
$\rho$ can be taken as the local form  (\ref{y2}) by Lemma
\ref{lem22} (ii). Therefore, to prove Theorem \ref{th2}, we only
need to solve $V$ from (\ref{y38})  with $h$ and $\rho$ being
given by (\ref{y2}).

Since $V$ is a conformal vector field of $h$ by the first equation
of (\ref{y38}), it follows from Lemma \ref{lem22} (i) that $V$ can
be locally expressed as
 \be\label{y45}
 V^i=-2\big(\tau+\langle
 \eta,x\rangle\big)x^i+|x|^2\eta^i+q_r^ix^r+\gamma^i,
 \ee
where $\tau$ is a constant number, $\eta,\gamma$ are constant
vectors and
  $Q=(q_i^j)$ is a skew-symmetric constant matrix. In this case, it
  follows from (\ref{g15}) that $c$ in (\ref{y38}) is given by
  \be\label{y46}
 c=\frac{\tau(1-\mu |x|^2)+ \langle \mu \gamma+\eta,x\rangle}{1+\mu |x|^2},
  \ee
The second equation in (\ref{y38}) is equivalent to
 \be\label{y47}
 V^j\frac{\pa p_i}{\pa
 x^j}+p_j\frac{\pa V^j}{\pa x^i}=-  2cp_i,
 \ee
and $p_i$ in (\ref{y2}) is given by
 \be\label{y48}
 p_i=\frac{4}{(1+\mu|x|^2)^2}\Big\{-2(\lambda+\mu
 \langle e,x\rangle) x^i+(1+\mu|x|^2)
 e^i\Big\}
 \ee
Now plugging (\ref{y45}), (\ref{y46}) and (\ref{y48}) into
(\ref{y47}) yields an equivalent equation of a polynomial in
$(x^i)$ of order four, in which every order must be zero. Then we
respectively have (from order zero to order four)
 \be\label{y49}
 Qe=-2\lambda \gamma,
 \ee
 \be\label{y50}
 \big(\langle
 \eta-\mu\gamma,e\rangle+2\lambda\tau\big)x^i-(\eta^i+\mu\gamma^i)\langle
 e,x\rangle=0,
 \ee
 \be\label{y51}
 \big[\lambda(\eta^i-\mu\gamma^i)-\mu
 q_k^ie^k\big]|x|^2+\mu\big[2\lambda\langle\gamma,x\rangle+4\tau\langle
 e,x\rangle-\langle e,Qx\rangle\big]x^i=0,
 \ee
 \be\label{y52}
 \big[\mu(\langle \eta-\mu\gamma, e\rangle
 -2\lambda\tau)x^i-\mu\langle
 e,x\rangle(\eta^i+\mu\gamma^i)\big]|x|^2+2\mu\langle
 e,x\rangle\langle\eta+\mu\gamma,x\rangle x^i=0,
 \ee
 \be\label{y53}
 \mu\big[(2\lambda \eta^i-\mu q^i_ke^k)|x|^2-2\langle
 2\lambda\eta-\mu Qe,x\rangle x^i\big]=0.
 \ee

\noindent{\bf Step A:} Assume $\mu=0$. Then by (\ref{y51}) we have
$\lambda=0$ or $\eta=0$.
 Contracting (\ref{y50}) by $x^i$ and using $\mu=0$, we obtain
  $$\big(\langle\eta,e\rangle+2\lambda\tau\big)|x|^2-\langle\eta,x\rangle\langle
 e,x\rangle=0,
 $$
 which implies
 \be\label{yg60}
\langle\eta,e\rangle+2\lambda\tau=0,\ \ \eta=0 \ or \ e=0.
 \ee
By (\ref{y49}) and (\ref{yg60}), it is easy to see that
(\ref{y49})--(\ref{y53}) are equivalent to one of the following
three cases:
 \be\label{y54}
 \lambda=0,\ \ \ \eta=0,\ \ \ Qe=0,\ \ \ (\mu=0);
 \ee
 \be\label{y55}
 \lambda=0,\ \ \  e=0,\ \ \ (\mu=0);
 \ee
 \be\label{y56}
 \eta=0,\ \ \  \tau=0,\ \ \ Qe=-2\lambda \gamma, \ \ \ (\mu=0).
 \ee

\noindent{\bf Step B:} Assume $\mu\ne0$. Then by (\ref{y53}) we
get
 \be\label{y57}
 2\lambda\eta-\mu Qe=0.
 \ee
Now plugging (\ref{y49}) into (\ref{y57}) yields
 \be\label{y58}
 \lambda=0,\ \ \ {\text or } \ \ \eta=-\mu \gamma.
 \ee
Contracting (\ref{y50}) by $x^i$ we have
 \be\label{y59}
 \langle
 \eta-\mu\gamma,e\rangle+2\lambda\tau=0,\ \ \ (\eta^i+\mu\gamma^i)
 e=0.
 \ee

\noindent{\bf Case I:} Assume $\lambda=0$ (see (\ref{y58})). Then
by (\ref{y59}) we have
 \be\label{y60}
 e=0,\ \ \ {\text or} \ \ \ \eta=-\mu \gamma\ \  {\text and} \ \
 \langle\eta,e\rangle= \langle\gamma,e\rangle=0.
 \ee

{\bf Case Ia:} Assume $e=0$ (see (\ref{y60})). In  this case, we
have
 \be\label{y61}
\lambda=0,\ \ \ e=0,\ \ (\mu\ne0),
 \ee
which implies (\ref{y49})--(\ref{y53}) automatically hold.

\

{\bf Case Ib:} Assume $\eta=-\mu \gamma$  and
 $\langle\eta,e\rangle= \langle\gamma,e\rangle=0$ ($e\ne 0$) (see (\ref{y60})). Plug $\lambda=0$
 and $Qe=0$ (see (\ref{y49})) into (\ref{y51}) and then we get
 $\tau=0$. So in this case, we have
  \be\label{y62}
\lambda=\tau=0,\ \ \ \langle\eta,e\rangle=
\langle\gamma,e\rangle=0,\ \ \ Qe=0,\ \ \ \eta=-\mu \gamma,\ \
(\mu\ne0),
  \ee
which implies (\ref{y49})--(\ref{y53}) automatically hold.

\

\noindent{\bf Case II:} Assume $\eta=-\mu \gamma$ (see
(\ref{y58})). By (\ref{y59}) we have
 \be\label{y63}
 \langle\gamma,e\rangle=\frac{\lambda\tau}{\mu}.
 \ee
Plugging (\ref{y49}) and $\eta=-\mu \gamma$ into (\ref{y51})
yields
 \be\label{y64}
 \tau=0,\ \ \ {\text or}\ \ \ e=0.
 \ee
Now it follows from $\eta=-\mu \gamma$, (\ref{y49}), (\ref{y63})
and (\ref{y64}) that
 \be\label{y65}
\eta=-\mu \gamma,\ \ \ \tau=\langle\gamma,e\rangle=0,\ \ \
Qe=-2\lambda\gamma,\ \ (\mu\ne0),
 \ee
which implies (\ref{y49})--(\ref{y53}) automatically hold.

Summing up the above discussions, we have six conditions:
(\ref{y54})--(\ref{y56}), (\ref{y61}), (\ref{y62}) and
(\ref{y65}). However, it is easy to show  (\ref{y62})
$\Rightarrow$ (\ref{y65}), and either (\ref{y55}) or (\ref{y61})
implies $F=\alpha$ is Riemannian since $\beta=\rho=0$ by
(\ref{y37}) and (\ref{y48}). Therefore, we have three independent
cases: (\ref{y54}), (\ref{y56}) and (\ref{y65}).   Now plug
(\ref{y54}), (\ref{y56}) and (\ref{y65}) into (\ref{y45})
respectively, we obtain (\ref{y3}), (\ref{y03}) and (\ref{y4}).
Then by (\ref{y46}), we can get $c$ in the three cases.

\

\noindent{\bf Condition (ii):} Suppose $\beta$ is  parallel with
respect to $\alpha$.

Since $F$ is locally projectively flat, $\alpha$ is also locally
projectively flat and then is of constant sectional curvature
$\mu$. If $\beta=0$, then $F=\alpha$ is Riemannian, which is
excluded. If $\beta\ne 0$, then $\mu=0$ and so $F$ is
flat-parallel ($\alpha$ is flat and $\beta$ is parallel). In this
case, $\beta$ is also a closed 1-form which is conformal with
respect to $\alpha$. By Proposition \ref{th1}, we only need to
solve $V$ from (\ref{y1}), which is just the case (\ref{y54}),
since we can put $b_i$ in the form of the right hand side of
(\ref{y48}) with $\mu=\lambda=0$, and thus $V$ is given by
(\ref{y3}).    \qed

\section{Homothetic fields of Randers or square metrics}\label{sec5}

In Theorem \ref{th2}, if $F=\alpha+\beta$ is a  Randers metric, we
can obtain the same result by further assuming that $F$ is of
isotropic S-curvature.

\begin{cor}\label{cor41}
 Let $F=\alpha+\beta$ be a locally projectively flat Randers
 metric of isotropic S-curvature. Then any conformal vector field $V$ of $F$
is homothetic, and $V$ can be locally written in the
 form (\ref{y3}), (\ref{y03}), or (\ref{y4}).
\end{cor}

{\it Proof :} Since $\beta$ is closed,  $F$ is of isotropic
S-curvature ${\bf S}=(n+1)\tau F$ for  some scalar function
$\tau=\tau(x)$ if and only if (\cite{CS})
 \be\label{y67}
 b_{i;j}=2\tau (a_{ij}-b_ib_j),
 \ee
where the covariant derivative is taken
  with respect to the Levi-Civita connection of $\alpha$. Define  a Riemann metric
$h=\sqrt{h_{ij}y^iy^j}$ and a 1-form $\rho=p_iy^i$ by
  \be\label{y66}
h:=\alpha,\ \ \rho:=\frac{1}{\sqrt{1-b^2}}\beta,\ \
 (b^2:=||\beta||_{\alpha}^2).
  \ee
 Under
  the deformation (\ref{y66}), it can be directly verified that
  (\ref{y67}) is equivalent to
   \be\label{y077}
 p_{i|j}=\frac{2\tau}{\sqrt{1-b^2}}\ h_{ij} \ \ (=-2ch_{ij}),
   \ee
where the covariant derivative is taken
  with respect to the Levi-Civita connection of $h$.
 It follows from (\ref{y077}) that $\rho$ is a closed
 1-form which is conformal with respect to $h$. Besides, $h$ is of
 constant sectional curvature (\cite{BaMa}). Thus we can express $h$ and $\rho$
locally by (\ref{y2}). Now by the proof shown in Theorem
\ref{th2}, we get the conclusions in this corollary.   \qed

\

Now we consider conformal vector fields of
$(\alpha,\beta)$-metrics which are locally projectively flat with
constant flag curvature. By \cite{LS} \cite{Y03}, if $F$ is an
$n(\ge 2)$-dimensional locally projectively flat
$(\alpha,\beta)$-metric with constant flag curvature, then $F$ is
flat-parallel, or $F=\alpha+\beta$ is a metric of Randers type, or
$F=(\alpha+\beta)^2/\alpha$ is metric of  square type.

\begin{lem}\label{lem42}
 Let $F$ be an $n(\ge 2)$-dimensional locally projectively flat Randers metric or  square
 metric. Define $h$ and $\rho$ by (\ref{y66}) for a Randers metric
 $F=\alpha+\beta$, or by $h:=(1-b^2)\alpha$ and
 $\rho:=\sqrt{1-b^2}\beta$ for a square metric
 $F=(\alpha+\beta)^2/\alpha$. Then $F$ is of constant flag
 curvature iff. $h$ and $\rho$ can be  locally given by
(\ref{y2}) with the relation
 \be\label{yg61}
 (1+4|e|^2)\mu+4\lambda^2=0.
 \ee
\end{lem}

We omit the proof of Lemma \ref{lem42}. By Theorem \ref{th2} and
Lemma \ref{lem42}, we can easily prove the following corollary.

\begin{cor}\label{cor42}
 Let $F$ be an $n(\ge 2)$-dimensional locally projectively flat Randers metric or  square
 metric with constant flag curvature.  Then any conformal vector field $V$
of $F$ is homothetic. Further, put $V=V^i(x)\pa/\pa x^i$, and then
related to $h$ and $\rho$ defined by (\ref{y2}) satisfying
(\ref{yg61}),
 locally  we have one of the following  cases:
 \ben
  \item[{\rm (i)}] ($\mu=0$, $\lambda=0$) $V$ is given by
   $$
 V^i=-2\tau x^i+q^i_kx^k+\gamma^i,\ \ \ (Qe=0),
   $$
  In this case,   $F$ is
  flat-parallel and
  $c$ in (\ref{y1}) is given by $c=\tau$.

 \item[{\rm (ii)}] ($\mu< 0$)  $V$ is given by
   $$
 V^i=2\mu\langle \gamma,x\rangle
 x^i+(1-\mu|x|^2)\gamma^i+q^i_kx^k,\ \ \ \big(Qe=-2\lambda \gamma\big),
   $$
  In this case,
  $c$ in (\ref{y1}) is given by $c=0$, and so $V$ is a Killing
  field.
 \een
 In the above,  $\tau$ is a constant number, $\gamma=(\gamma^i)$ is a constant
 vector, and $Q=(q^i_k)$ is a constant skew-symmetric matrix.
\end{cor}

\vspace{0.6cm}

\noindent Guojun Yang \\
Department of Mathematics \\
Sichuan University \\
Chengdu 610064, P. R. China \\
yangguojun@scu.edu.cn

\end{document}